\documentclass[12pt,a4paper]{paper}
\usepackage{amsmath}
\usepackage{amssymb}
\usepackage{amsthm}
\usepackage{amsfonts}
\usepackage{latexsym}

\theoremstyle{plain}
\newtheorem{thm}{Theorem}
\newtheorem{cor}{Corollary}[section]
\newtheorem{lem}{Lemma}[section]
\newtheorem{prop}{Proposition}[section]
\theoremstyle{definition}
\newtheorem{defn}{Definition}[section]

\theoremstyle{remark}

\newtheorem{rem}{Remark}[section]

\title{Semiclassical almost isometry}
\author{Roberto Paoletti
}
\date{}

\begin{document}

\maketitle

\noindent \footnotesize{{\bf Address.} Dipartimento di Matematica e
Applicazioni, Universit\`{a} di Milano Bicocca, Via R. Cozzi 53,
20126 Milano, Italy; {\bf e-mail:} roberto.paoletti@unimib.it}

\abstract{Let $(M,\omega,J)$ be a compact and connected polarized
Hodge manifold, $\mathfrak{S}$ an isodrastic leaf of half-weighted
Bohr-Sommerfeld Lagrangian submanifolds. We study the relation
between the Weinstein symplectic structure of $\mathfrak{S}$ and the
asymptotics of the the pull-back of the Fubini-Study form under the
projectivization of the so-called BPU maps on $\mathfrak{S}$.}



\section{Introduction}

Let $(M,J)$ be an irreducible n-dimensional complex projective
manifold, $A\rightarrow M$ an ample line bundle, $h$ an Hermitian
metric on $A$ such that the curvature of the unique compatible
covariant derivative is $-2\pi i \omega$, where $\omega$ is a
K\"{a}hler form on $M$. By the Tian-Zelditch almost isometry theorem
\cite{tian}, \cite{z}, the projective embeddings $\varphi
_k=:\varphi_{A^{\otimes k}}:M\rightarrow \mathbb{P}\big
(H^0(M,A^{\otimes k})^*\big )$ are asymptotically symplectic as
$k\rightarrow +\infty$, in an appropriate rescaled sense. Thus the
symplectic structure of the classical phase space $(M,\omega)$ is
encapsulated in the asymptotics of its quantizations
$H^0(M,A^{\otimes k})$. However, in light of the uncertainty
principle and of the WKB method, the geometric objects most
naturally associated to quantum physical states (the \textit{true
points} of phase space) are Lagrangian submanifolds of $(M,\omega)$
\cite{bw}, \cite{gs-gq}, \cite{gs-gc}, \cite{weinst-hb}. Motivating
the study of quantization and reduction, this point of view led to
Weinstein's discovery of a natural symplectic structure on
isodrastic leaves of weighted Lagrangian submanifolds
\cite{weinst-hb}. One may then ask whether almost isometry still
holds when $(M,\omega)$ is replaced by an isodrastic leaf of compact
and connected half-weighted Bohr-Sommerfeld Lagrangian submanifolds,
endowed with a closed 2-form $\Omega$ of Weinstein type, and the
$\varphi _k$'s by their semiclassical analogues taking value in
$\mathbb{P}H^0(M,A^{\otimes k})\cong \mathbb{P}^{N_k}$, and denoted
$\Phi_k$ below; these are the (projectivisation of the) maps
introduced in \cite{bpu} (and called BPU maps in \cite{gt}).


Let $\mathfrak{S}$ be the manifold of all half-weighted
Bohr-Sommerfeld Lagrangian submanifolds $(L,\lambda)$ of
$(M,\omega)$ such that $L$ is isodrastically equivalent to a given
$L_0$ (\S \ref{subsctn:weighted-lag-sub}, \S
\ref{subsectn:half-weighted}). Besides $\Omega$, the isodrastic leaf
$\mathfrak{S}$ also carries a natural semidefinite Riemannian metric
$G$. In fact, $G$ and $\Omega$ are non-degenerate and compatible,
hence define an almost K\"{a}hler structure, on the open subset
$\mathfrak{S}'\subseteq \mathfrak{S}$ of all pairs $(L,\lambda)\in
\mathfrak{S}$ with $\lambda$ nowhere vanishing on $L$. Given the
K\"{a}hler structure $(M,J,\omega)$, the smooth tangent space
$T_{(L,\lambda)}\mathfrak{S}$ to $\mathfrak{S}$ at $(L,\lambda)$ is
naturally isomorphic to the space of pairs $(f,\ell)$, where $f\in
\mathcal{C}^\infty(L)$ and $\ell$ is a $\mathcal{C}^\infty$
half-density on $L$, satisfying $\int _Lf\lambda\bullet \lambda=\int
_L\ell \bullet \lambda=0$ (statement i) of Proposition
\ref{prop:tg-space-half-weights}). Let $W(f,\ell)$ be the tangent
vector associated to a pair $(f,\ell)$. For every $k\gg 0$ suitably
divisible, let $\Phi _k:U_k\rightarrow \mathbb{P}^{N_k}$ be the
$k$-th projectivized BPU map (\S \ref{subsctn:proj-bpu-maps});
$U_k\subseteq \mathfrak{S}$ is an appropriate open subset, and
$(L,\lambda)\in U_k$ for any $(L,\lambda)\in \mathfrak{S}$ and for
all suitably divisible $k\gg 0$. Thus, $\mathfrak{S}=\bigcup _kU_k$.

In Weinstein's setting, almost isometry does not hold literally (if
anything because an isodrastic leaf is infinite dimensional).
Nonetheless, $\Omega$ and $G$ can be extracted from the asymptotics
of BPU maps:

\begin{thm}\label{thm:semicl-asympt-isom}
Let $\mathfrak{S}$ be an isodrastic leaf of half-weighted compact
and connected Bohr-Sommerfeld Lagrangian submanifolds of
$(M,\omega)$. For every $(L,\lambda)\in \mathfrak{S}$ and
$W(f,\ell),W(f',\ell')\in T_{(L,\lambda)}\mathfrak{S}$, the
following asymptotic expansions hold as $k=lr$ and $l\rightarrow
+\infty$:
\begin{eqnarray*}
\Phi_k^*(\omega_{FS})_{(L,\lambda)}\Big(W(f,k\ell),W(f',k\ell')\Big)&\sim&k^2\cdot
\Omega _{(L,\lambda)}\Big (W(f,l),W(f',\ell')\Big)\\
&&+\sum _{h\ge 1}b_h(L,\lambda,f,\ell,f',\ell')\,k^{2-h/2},
\end{eqnarray*}
\begin{eqnarray*}
\Phi_k^*(g_{FS})_{(L,\lambda)}\Big(W(f,k\ell),W(f',k\ell')\Big)&\sim&k^2\cdot
G _{(L,\lambda)}\Big (W(f,l),W(f',\ell')\Big)\\
&&+\sum _{h\ge 1}c_h(L,\lambda,f,\ell,f',\ell')\,k^{2-h/2}.
\end{eqnarray*}
Here, $r\in \mathbb{N}$ is the an invariant of $\mathfrak{S}$ given
by order of the image in $S^1$ of the holonomy representation of $L$
for the unit circle bundle $X\subseteq A^*$. If $s\ge 1$,
$S\mathfrak{S}\subseteq T\mathfrak{S}$ is the unit sphere bundle
(for any given smooth metric), and $K\subseteq S\mathfrak{S}$ is the
image of a smooth map from a compact subset of $\mathbb{R}^s$, these
asymptotic expansions are uniform on the set of pairs of tangent
vectors multiples of elements of $K$.
\end{thm}

We hint at two possible extensions: to the category of compact
almost K\"{a}hler manifolds, on the hand, in view of the microlocal
description of almost complex Szeg\"{o} kernels in \cite{sz}, and on
the other allowing $L$ to be open, but requiring the half-weights to
be compactly supported. We shall leave these generalizations to the
interested reader.

\bigskip

\noindent \textbf{Ackowledgments.} I am indebted to the the referee
for suggesting some improvements in presentation.

\section{Preliminaries}

We shall denote by $\mathcal{D}^\infty(Z)$ (resp.,
$\mathcal{D}^\infty_{(1/2)}(Z)$) the space of all
$\mathcal{C}^\infty$ real-valued densities (resp., half-densities )
on a manifold $Z$.
There is a natural commutative product
$\bullet:\mathcal{D}_{(1/2)}^\infty(Z)\otimes
\mathcal{D}_{(1/2)}^\infty(Z)\rightarrow \mathcal{D}^\infty(Z)$,
$\lambda \otimes \eta\mapsto \lambda \bullet \eta$, given by
pointwise multiplication of functions on frame bundles. All
densities and half-densities will be understood to be real-valued.
Given a Riemannian structure on $Z$, $\mathrm{dens}_Z$ (resp.,
$\mathrm{dens}_Z^{(1/2)}$) will denote the corresponding volume
density (resp., half-density).

\subsection{Weighted Lagrangian
and Planckian submanifolds}\label{subsctn:weighted-lag-sub}

The space $\mathbf{L}=\mathbf{L}(M,\omega)$ of all compact and
connected Lagrangian submanifolds of $(M,\omega)$ is an
infinite-dimensional manifold. The (smooth) tangent space of
$\mathbf{L}$ at any $L\in \mathbf{L}$ is
$T_L\mathbf{L}=\mathcal{Z}^1(L)$, the vector space of all closed
1-forms on $L$. Furthermore, $\mathbf{L}$ carries a natural
integrable distribution $\mathfrak{B}\subseteq T\mathbf{L}$, whose
value at any $L\in \mathbf{L}$ is the subspace
$\mathcal{B}^1(L)\subseteq \mathcal{Z}^1(L)$ of all exact 1-forms on
$L$. $\mathfrak{B}$ is called the \textit{isodrastic} distribution,
and its leaves the isodrastic leaves of $\mathbf{L}$. Lagrangian
submanifolds $L,L'\in \mathbf{L}$ belong to the same isodrastic leaf
(in which case they are called isodrastically equivalent) if and
only $L'$ can be deformed into $L$ by flowing it along globally
defined Hamiltonian vector fields  \cite{weinst-hb}. A compact and
connected \textit{weighted} Lagrangian submanifold of $(M,\omega)$
is a pair $(L,\varrho)$, where $L\in \mathbf{L}$ and $\varrho \in
\mathcal{D}^\infty (L)$ is a \textit{weight} on $L$, that is, $\int
_L\varrho =1$. We shall denote by
$\mathbf{WL}=\mathbf{WL}(M,\omega)$ the manifold of all such pairs.
Given the natural forgetful projection, $p:\mathbf{WL}\rightarrow
\mathbf{L}$, for any isodrastic leaf $\mathfrak{I}\subseteq
\mathbf{L}$ set $\mathbf{W}\mathfrak{I}=:p^{-1}\big
(\mathfrak{I}\big )$ (really an immersed submanifold). It is the
infinite dimensional manifold $\mathbf{W}\mathfrak{I}$ that carries
a built-in symplectic structure $\Omega _{\mathrm{Wein}}$ (\S 3 of
\cite{weinst-hb}).


\begin{defn}\label{defn:planckian}
Let $X\subseteq A^*$ be the unit circle bundle, with projection $\pi
:X\rightarrow M$, so that the connection 1-form $\alpha$ on $X$ is a
contact structure, and $d\alpha =\pi^*(\omega)$. A submanifold
$P\subseteq X$ is \textit{Planckian} if it is Legendrian and
furthermore (by restriction of $\pi$) an unramified cover of a
Lagrangian submanifold $L=\pi (P)\subseteq M$.
$\mathbf{P}=\mathbf{P}(X,\alpha)$ will denote the collection of all
compact and connected Planckian submanifolds of $(X,\alpha)$.
\end{defn}


\begin{defn}\label{defn:bohr-sommerfeld}
A submanifold $L\subseteq M$ is a \textit{Bohr-Sommerfeld Lagrangian
submanifold} (BSL for short) if $L=\pi (P)$ for some Planckian
submanifold $P\subseteq X$. Let $\mathbf{L}_{\mathrm{BS}}\subseteq
\mathbf{L}$ be the subspace of all compact and connected BSL
submanifolds.
\end{defn}


\begin{rem} \label{rem:holonomy-planckian}
Suppose $L\in \mathbf{L}$; then $L\in \mathbf{L}_{\mathrm{BS}}$ if
and only if the image of the holonomy representation $\pi
_1(L)\rightarrow S^1$ for the principal $S^1$-bundle $X$ is a finite
subgroup $\mathrm{Hol}(L)\subseteq S^1$. If $P\in \mathbf{P}$ is
such that $L=\pi(P)$, the projection $P\rightarrow L$ is an
unramified cover of degree $r=:\left |\mathrm{Hol}(L)\right
|$.\end{rem}


The property of being BSL is invariant under isodrastic deformations
\cite{weinst-hb}, hence $\mathbf{L}_{\mathrm{BS}}$ is a union of
isodrastic leaves. Define $\widetilde{\pi}:\mathbf{P}\rightarrow
\mathbf{L}_{\mathrm{BS}}$ by $\widetilde{\pi}(P)=:\pi(P)$. For any
isodrastic leaf $\mathfrak{I}\subseteq \mathbf{L}_{\mathrm{BS}}$,
$\mathbf{P}_{\mathfrak{I}}=:\widetilde{\pi}^{-1}\big
(\mathfrak{I}\big )$ is an infinite-dimensional manifold, and
$T_P\big (\mathbf{P}_{\mathfrak{I}}\big )\cong \mathcal{C}^\infty
(L)$ for any $P\in \mathbf{P}_{\mathfrak{I}}$, where $L=\pi (P)$
(\cite{weinst-hb}, Lemma 4.1). Furthermore, the image of the
holonomy representation $\pi_1(L)\rightarrow S^1$ associated to the
principal $S^1$-bundle $X$ is the same $\forall\,L\in \mathfrak{I}$;
its cardinality equals the degree of the unramified cover
$P\rightarrow L=:\pi (P)$, $\forall\,P\in \mathfrak{I}$. Denote this
image by $\mathrm{Hol}(\mathfrak{I})\subseteq S^1$, and set
$G_\mathfrak{I}=:S^1/\mathrm{Hol}(\mathfrak{I})\cong S^1$. Thus
$\mathbf{WP}_\mathfrak{I}=:\mathbf{P}_{\mathfrak{I}}\times
_{\mathfrak{I}}\mathbf{W}\mathfrak{I}$ consists of all pairs
$(P,\varrho)$, where $P\in \mathbf{P}_\mathfrak{I}$ and $\varrho$ is
a weight on $\pi (P)\in \mathfrak{I}$. The projection
$\widehat{\pi}:\mathbf{WP}_\mathfrak{I}\rightarrow
\mathbf{W}\mathfrak{I}$, given by $(P,\varrho)\mapsto
\big(\widetilde{\pi}(P),\varrho\big )$, is a principal
$G_\mathfrak{I}$-bundle, and has an intrinsic  \textit{universal}
connection in the terminology of \cite{weinst-hb}; the normalized
curvature of this connection is the symplectic structure
$\Omega_{\mathrm{Wein}}$ on $\mathbf{W}\mathfrak{I}$ (Proposition
4.3 of \cite{weinst-hb}). Heuristically, the circle bundle
$\widehat{\pi}:\mathbf{WP}_\mathfrak{I}\rightarrow
\mathbf{W}\mathfrak{I}$ is a semiclassical analogue of the circle
bundle $\pi :X\rightarrow M$. In the present K\"{a}hler context the
theory of \cite{weinst-hb} implies that the tangent space
$T_{(L,\varrho)}\mathbf{W}\mathfrak{I}$ has a simple intrinsic
description. Pairing the proof of Theorem 3.32 of \cite{ms} with
that of Lemma 3.14 of \cite{ms} yields the following:

\begin{lem}
\label{lem:normal-cotangent-neighborhood} Let $L\subseteq M$ be a
Lagrangian submanifold, $T^*L$ its cotangent bundle, with projection
$q:T^*L\rightarrow L$, and canonical symplectic structure $\omega
_{\mathrm{can}}$. Let $\mathbb{O}(L)=:\{(l,0):l\in L\}\subseteq
T^*L$. Then there exist open neighborhoods $L\subseteq U\subseteq M$
and $\mathbb{O}(L)\subseteq V\subseteq T^*L$, and a natural choice
of a symplectomorphism $\gamma:U\rightarrow V$, such that $\gamma
(l)=l$, and the inverse image $\gamma
^{-1}\Big(q^{-1}(l)\Big)\subseteq U$ of $q^{-1}(l)\subseteq T^*L$ is
perpendicular to $L$ at $l$ (for the K\"{a}hler metric), $\forall\,
l\in L$.\end{lem}

We shall call $\gamma$ the \textit{normal cotangent structure} of
$M$ near $L$, and denote the projection by $\beta:U\rightarrow L$.
The discussion surrounding equation (3) of \cite{weinst-hb} implies:

\begin{prop}
Given the (almost) K\"{a}hler structure $(M,\omega,J)$,

i): $\forall\,(L,\varrho)\in \mathbf{W}\mathfrak{I}$,
$T_{(L,\varrho)}\big (\mathbf{W}\mathfrak{I})$ is naturally
isomorphic to the vector space of all pairs $(f,\phi)\in
\mathcal{C}^\infty(L)\times \mathcal{D}^\infty(L)$ satisfying $\int
_Lf\,\varrho =\int _L\phi=0$;

ii): $\forall\,(P,\varrho)\in \mathbf{WP}_\mathfrak{I}$,
$T_{(P,\varrho)}\big (\mathbf{WP}_\mathfrak{I})$ is naturally
isomorphic to the vector space of all pairs $(f,\phi)\in
\mathcal{C}^\infty(L)\times \mathcal{D}^\infty(L)$ satisfying $\int
_L\phi=0$. \label{prop:tg-space-weights}            \end{prop}

In i), $f$ is implicitly extended by pull-back to $U$ under $\beta$,
and thus defines a Hamiltonian vector field $\upsilon _f$ on $U$.
The flow of $\upsilon_f$ determines a path $L_t$, $|t|<\epsilon$, of
Lagrangian submanifolds with $L_0=L$. The condition on $f$ is a
renormalization that fixes it uniquely. Restricting $\beta$
determines diffeomorphisms $\beta_t:L_t\rightarrow L$. If $\rho _t$
is a family of weights on $L$ with $\rho _0=\varrho$, we obtain by
pull-back weights $\varrho _t=\beta_t^*(\rho_t)$ on $L_t$. The first
order datum at $t=0$ corresponding to the family of weights
$\varrho_t$ is the density $\phi=\rho'(0)$ on $L$, so that $\rho
_t=\varrho +t\cdot \phi +O(t^2)$; the condition on $\phi$ results
from differentiating the constraint $\int_{L_t}\varrho _t=\int
_L\rho _t=1$ at $t=0$. If $(L,\varrho)\in \mathbf{W}\mathfrak{I}$
and $V_i=V(f_i,\phi _i)\in T_{(L,\varrho)}\big
(\mathbf{W}\mathfrak{I})$, $i=1,2$, by (3) on page 139 of
\cite{weinst-hb} their Weinstein symplectic pairing is
\begin{equation}\label{eqn:symp-pairng}
(\Omega _{\mathrm{Wein}})_{(L,\varrho)}\big (V_1,V_2\big )=\int
_L\big (f_1\,\phi _2-f_2\,\phi _1\big).\end{equation}

\begin{defn}\label{derivative-of-density-along-flow}
Suppose $L\in \mathbf{L}$, $f\in \mathcal{C}^\infty(L)$. For
sufficiently small $t$, let $\phi _t$ be the flow of $\upsilon _f$,
defined locally near $L$, and $L_t=:\phi _t(L)$.  Let
$\mathrm{dens}_{L_t}^{(1/2)}$ be the Riemannian half-density on
$L_t$, $\psi _t:L\rightarrow L_t$ the diffeomorphism $m\mapsto \phi
_t(m)$, and $\gamma
_t=:\psi_t^*\big(\mathrm{dens}_{L_t}^{(1/2)}\big)$. Thus $\gamma
_t=G_t\cdot\mathrm{dens}_{L}^{(1/2)}$ for a unique
$G_t\in\mathcal{C}^\infty(L)$. Then $\Gamma (L,f)=:\left
.\frac{\partial G_t}{\partial t}\right |_{t=0}\in
\mathcal{C}^\infty(L)$ depends linearly on $f$.
\end{defn}


\subsection{Half-weighted Lagrangian
submanifolds}\label{subsectn:half-weighted}

If $L\in \mathbf{L}$, $\lambda \in \mathcal{D}_{(1/2)}^\infty(L)$ is
a \textit{half-weight} if $\lambda \bullet \lambda $ is a weight on
$L$. Let
$\mathbf{W}_{\mathrm{h}}\mathbf{L}=\mathbf{W}_{\mathrm{h}}\mathbf{L}(M,\omega)$
be the space of all compact and connected half-weighted Lagrangian
submanifolds of $(M,\omega)$, and by
$\mathbf{W}_{\mathrm{h}}\mathbf{P}=\mathbf{W}_{\mathrm{h}}\mathbf{P}(X,\alpha)$
the space of all pairs $(P,\lambda)$, where $P\in \mathbf{P}$ and
$\lambda$ is a half-weight on $\pi (P)$. Given an isodrastic leaf
$\mathfrak{I}\subseteq \mathbf{L}_{\mathrm{BS}}$, we have leaves
$\mathbf{W}_{\mathrm{h}}\mathfrak{I}$ and
$\mathbf{W}_{\mathrm{h}}\mathbf{P}_{\mathfrak{I}}$. The analogue of
Proposition \ref{prop:tg-space-weights} is:

\begin{prop}\label{prop:tg-space-half-weights}
Given the (almost) K\"{a}hler structure $(M,\omega,J)$,

i): $\forall\,(L,\lambda)\in \mathbf{W}_{\mathrm{h}}\mathfrak{I}$,
$T_{(L,\lambda)}\big (\mathbf{W}_{\mathrm{h}}\mathfrak{I}\big)$ is
naturally isomorphic to the vector space of all pairs $(f,\ell)\in
\mathcal{C}^\infty(L)\times \mathcal{D}^\infty(L)$ satisfying $\int
_Lf\,\lambda \bullet \lambda =\int _L\ell \bullet \lambda=0$;

ii): $\forall\,(P,\lambda)\in
\mathbf{W}_{\mathrm{h}}\mathbf{P}_{\mathfrak{I}}$,
$T_{(P,\lambda)}\big
(\mathbf{W}_{\mathrm{h}}\mathbf{P}_{\mathfrak{I}})$ is naturally
isomorphic to the vector space of all pairs $(f,\ell)\in
\mathcal{C}^\infty(L)\times \mathcal{D}^\infty(L)$ satisfying $\int
_L\ell \bullet \lambda = 0$. \end{prop}

\begin{defn} If $(f,\ell)\in
\mathcal{C}^\infty(L)\times \mathcal{D}^\infty(L)$ satisfy $\int
_Lf\,\lambda \bullet \lambda =\int _L\ell \bullet \lambda=0$, let
$W(f,\ell)\in T_{(L,\lambda)}\big
(\mathbf{W}_{\mathrm{h}}\mathfrak{I}\big)$ be the tangent vector
associated to $(f,\ell)$. Similarly, if $(f,\ell)\in
\mathcal{C}^\infty(L)\times \mathcal{D}^\infty(L)$ satisfy $\int
_L\ell \bullet \lambda = 0$, let $\widetilde{W}(f,\ell)\in
T_{(P,\lambda)}\big
(\mathbf{W}_{\mathrm{h}}\mathbf{P}_{\mathfrak{I}})$ be the tangent
vector associated to $(f,\ell)$. If $P\in \mathbf{P}_{\mathfrak{I}}$
is Planckian and $L=\pi(P)$ any $W(f,\ell)\in T_{(L,\lambda)}\big
(\mathbf{W}_{\mathrm{h}}\mathfrak{I}\big)$ lifts to
$\widetilde{W}(f,\ell)\in T_{(P,\lambda)}\big
(\mathbf{W}_{\mathrm{h}}\mathbf{P}_{\mathfrak{I}})$.\end{defn}

Define $\Psi:\mathbf{W}_{\mathrm{h}}\mathfrak{I}\rightarrow
\mathbf{W}\mathfrak{I}$ by $(L,\lambda)\mapsto (L,\lambda \bullet
\lambda)$. The differential $d_{(L,\lambda)}\Psi:T_{(L,\lambda)}\big
(\mathbf{W}_{\mathrm{h}}\mathfrak{I}\big )\rightarrow
T_{(L,\lambda\bullet \lambda)}\big (\mathbf{W}\mathfrak{I}\big )$ is
given by $ d_{(L,\lambda)}\Psi \big ((f,\ell)\big)=(f,2\ell \bullet
\lambda)$. Clearly, $\mathbf{W}_{\mathrm{h}}\mathbf{P}_\mathfrak{I}$
is the principal $S^1$-bundle on
$\mathbf{W}_{\mathrm{h}}\mathfrak{I}$ obtained by pulling back
$\widehat{\pi}:\mathbf{WP}_{\mathfrak{I}}\rightarrow
\mathbf{W}{\mathfrak{I}}$ by $\Psi$. We shall also write
$\widehat{\pi}$ the projection
$\mathbf{W}_{\mathrm{h}}\mathbf{P}_\mathfrak{I}\rightarrow
\mathbf{W}_{\mathrm{h}}\mathfrak{I}$.  Consider the closed 2-form
$\Omega=:\Psi ^*\big(\Omega_{\mathrm{Wein}}\big)$ on
$\mathbf{W}_{\mathrm{h}}\mathfrak{I}$. If $(L,\lambda)\in
\mathbf{W}_{\mathrm{h}}\mathfrak{I}$, $W=W(f,\ell ),
W'=W(f',\ell')\in T_{(L,\lambda)}\big
(\mathbf{W}_{\mathrm{h}}\mathfrak{I}\big )$, we have:
\begin{equation}\label{eqn:symp-str-half}
\Omega_{(L,\lambda)}(W,W')=2\int _L\big (f\cdot\ell '-f'\cdot \ell
\big)\bullet \lambda.\end{equation}

\noindent Now on $\mathbf{W}_{\mathrm{h}}\mathfrak{I}$ we also have
a positive semi-definite Riemannian structure, given by
\begin{equation}
\label{eqn:riem-str-half} G_{(L,\lambda)}(W,W')=2\int _L\Big
((f\,f')\cdot \lambda \bullet \lambda +\ell\bullet \ell'\Big).
\end{equation}
Endowed with $\Omega$ and $G$, the open subset
$\mathbf{W}'_{\mathrm{h}}\mathfrak{I}\subseteq
\mathbf{W}_{\mathrm{h}}\mathfrak{I}$ of pairs $(L,\lambda)\in
\mathbf{W}_{\mathrm{h}}\mathfrak{I}$ such that $\lambda$ is nowhere
vanishing is an infinite-dimensional almost K\"{a}hler manifold. In
fact, $\Omega$ and $G$ are non-degenerate on
$\mathbf{W}'_{\mathrm{h}}\mathfrak{I}$, and related by the almost
complex structure
$\mathbb{J}\Big(W(f,g\,\lambda)\Big)=W(-g,f\,\lambda)$.

\subsection{Good coordinates along a Legendrian submanifold}

For $r\in \mathbb{N}$ and $\epsilon >0$, let $B^r_\epsilon\subseteq
\mathbb{R}^r$ be the open ball of radius $\epsilon$ centered at the
origin. Let $(p,q,\theta)$ be the standard linear coordinates on
$\mathbb{R}^{2\mathrm{n}+1}\cong \mathbb{R}^{\mathrm{n}}\times
\mathbb{R}^{\mathrm{n}}\times \mathbb{R}$; if $\psi
:B^{2\mathrm{n}+1}_\epsilon\rightarrow V=:\psi \big
(B^{2\mathrm{n}+1}_\epsilon\big)\subseteq X$ is a local chart,
$(p,q,\theta):V\rightarrow \mathbb{R}^{2\mathrm{n}+1}$ will also
denote the induced local coordinates, and $\frac{\partial}{\partial
q_i}, \,\frac{\partial}{\partial p_i},\,\frac{\partial}{\partial
\theta}$ the corresponding vector fields on $V$.

The compatible connection defines a direct sum decomposition
$TX=\mathrm{Hor}(X)\oplus \mathrm{Ver}(X)$, where
$\mathrm{Hor}(X)=\ker (\alpha)$, $\mathrm{Ver}(X)=\ker (d\pi)$. To
define a Riemannian structure $g_X$ on $X$,  we declare this to be
an orthogonal direct sum, take the pull-back of the Riemannian
structure on $M$ as a metric on $\mathrm{Hor}(X)$, and require the
generator of the $S^1$-action to have norm $\frac{1}{2\pi}$. The
$S^1$-orbits have unit length for $g_X$, and for the corresponding
Riemannian density $\mathrm{dens}_X$, the natural isomorphism
$H^0\big (M,A^{\otimes k}\big )\cong H(X)_k\subseteq L^2(X)$ is
unitary; here $H(X)_k$ is the level $k$ Hardy space of $X$. Given
$x\in X$, let $\|\cdot\|_x$ be associated the norm on $T_xX$.


\begin{prop}
\label{prop:local-adapted-coordinates} Let $P \subseteq X$ be a
Legendrian submanifold. For any $x\in P$, there exists a local chart
$\psi :B^{2\mathrm{n}+1}_\epsilon\rightarrow X$ for $X$, centered at
$x$ and such that:
\begin{enumerate}
  \item $P \cap V$ is defined by the conditions $p=0$ and $\theta=0$, where
  $V=:\psi \big
(B^{2\mathrm{n}+1}_\epsilon\big)$;
  \item $e^{i\vartheta}\cdot \psi (p,q,\theta)=\psi (p,q,\theta
  +\vartheta)$, whenever all terms are defined;
  \item at any $y\in V\cap P$,
  $\mathrm{Hor}(X/M)_y=\mathrm{span}\left \{\left .\frac{\partial}{\partial q_i}\right |_y, \,
  \left .\frac{\partial}{\partial p_i}\right |_y\,:\,1\le i\le\mathrm{n}\right
  \}$;
  \item for every $y\in V\cap P$, one has
  \begin{eqnarray*}\mathrm{span}\left \{\left .
  \frac{\partial}{\partial p_i}\right |_y\,:\,1\le i\le\mathrm{n}\right \}
  =J_y
  \left (\mathrm{span}\left \{\left .\frac{\partial}{\partial q_i}\right |_y\,:\,1\le i\le\mathrm{n}
  \right \}\right )=J_y
  \left (T_yP\right ),\end{eqnarray*}
  where $J_y\in \mathrm{End}\big(\mathrm{Hor}(X/M)_y\big)$ is induced by the complex
  structure of $M$;
  \item if $y\in V\cap P$ and $\eta_1,\ldots,\eta_\mathrm{n}\in
  \mathbb{R}$, then
  $ \left \|\sum _{j=1}^{\mathrm{n}} \eta_j\,\left.\frac{\partial}{\partial p_j}\right|_y\right
  \|_y^2=\sum _{j=1}^{\mathrm{n}}\eta _j^2.$
\end{enumerate}
\end{prop}


\textit{Proof.} In the following, $\epsilon>0$ is allowed to vary
from line to line. By \S 2 of \cite{dp}, for any $y\in P$ there
exists a system of Heisenberg local coordinates
$(p^{(y)},q^{(y)},\theta^{(y)})$ adapted to $P$ at $y$. This means
that $(p^{(y)},q^{(y)},\theta^{(y)})$ are local Heisenberg
coordinates for $X$ centered at $y$, in the sense of \cite{sz}, and
that $P$ is tangent to the locus $p^{(y)}=0$ at $y$. This
construction may be deformed smoothly with $y\in P$: $\forall\,x\in
P$ there exist an open neighborhood $x\in P'\subseteq P$, $\epsilon
>0$, and a smooth map $\Psi :P' \times
B^{2\mathrm{n}}_\epsilon\times (-\pi,\pi)\rightarrow X$, such that
$\forall\,y\in P'$ the partial map
$\Psi(y,\cdot,\cdot):B^{2\mathrm{n}}_\epsilon\times
(-\pi,\pi)\rightarrow X$ is an Heisenberg local chart adapted to $P$
at $y$. We may assume without loss that $P'$ is the image of a local
chart $\phi :B^{\mathrm{n}}_\epsilon\rightarrow P'$ for $P$ centered
at $x$. Let $q=(q_i)$ denote the linear coordinates on
$\mathbb{R}^\mathrm{n}$. Define $\psi
:B^{\mathrm{2n}+1}_\epsilon\rightarrow X$ by $ \psi
(p,q,\theta)=\Psi \big(\phi(q),(p,0,\theta)\big)$.
By definition of Heisenberg local coordinates, $\psi$ is a local
chart for $X$ satisfying all the conditions in the statement of the
Proposition. \hfill $\square$

\begin{defn}\label{defn:good-local-coordinates}
A system of local coordinates defined as in Proposition
\ref{prop:local-adapted-coordinates} will be called \textit{a system
of good local coordinates for $X$ along $P$}.
\end{defn}

In the notation of the Proposition, $\pi (V)\subseteq M$ is an open
subset, $\pi (P\cap V)\subseteq \pi(V)$ is a Lagrangian submanifold,
and $(p,q)$ is naturally a local coordinate chart on $\pi(V)$, in
which $\pi (P\cap V)$ is defined by the condition $p=0$. The
Heisenberg local charts $\Psi(y,\cdot,\cdot)$ appearing in the proof
are defined on $B^{2\mathrm{n}}_\epsilon\times (-\pi,\pi)$, but the
good coordinate chart is defined on
$B^{\mathrm{2n}+1}_{\epsilon}\subseteq
B^{2\mathrm{n}}_\epsilon\times (-\pi,\pi)$. This ensures that $P$
intersect each $S^1$-orbit at most once in the given chart, and the
image of $P\cap V$ in $X$ is a submanifold. Now suppose that
actually $P\in \mathbf{P}$. Let $r\in \mathbb{N}$ be the degree of
the unramified cover $P\rightarrow L=:\pi(P)\subseteq M$. Then
$e^{i\theta}\cdot P=P$ if $e^{i\theta}\in \mathbb{Z}_r=\left
<e^{2\pi i/r}\right
>\subseteq S^1$, and $\left (e^{i\theta}\cdot P\right )\cap
P=\emptyset$ if $e^{i\theta}\not\in \mathbb{Z}_r$. In fact,
$\mathbb{Z}_r$ acts as a group of Riemannian covering maps for
$P\rightarrow L=:\pi(P)$, and Proposition
\ref{prop:local-adapted-coordinates} may be strengthened:

\begin{prop}\label{prop:local-adapted-coordinates-planckian}
Suppose $P\in \mathbf{P}$. For any $x\in P$, there exists a local
chart $\psi :B^{2\mathrm{n}}_\epsilon\times (-\pi,\pi)\rightarrow X$
for $X$, such that $x=\psi (0,0,\vartheta_1)$ for some
$\vartheta_1\in (-\pi,\pi)$, satisfying conditions 2., 3., and 4. of
Proposition \ref{prop:local-adapted-coordinates} with $P$ in place
of $\Lambda$, and such that in addition condition 1 is replaced by:
\begin{description}
  \item[1a] let
  $V=:\psi \Big
(B^{2\mathrm{n}}_\epsilon\times (-\pi,\pi)\Big)$; then $P \cap V$ is
defined by the conditions $p=0$ and $\theta\in
\{\vartheta_1,\cdots,\vartheta_r\}$,
  where $\vartheta_j\in (-\pi,\pi)$ are all distinct;
  \item[1b] $P \cap
V=\pi ^{-1}\big(P \cap V\big)\cap \Lambda$.
\end{description}
\end{prop}

\subsection{Projectivized BPU maps}\label{subsctn:proj-bpu-maps}

Given the volume form $\mathrm{vol}_X=\alpha\wedge
\pi^*(\omega)^{\wedge \mathrm{n}}$, we shall identify functions,
densities and half-densities on $X$. Any $(P,\lambda)\in
\mathbf{W}_{\mathrm{h}}\mathbf{P}_\mathfrak{I}$ induces a
generalized half-density $\delta _{(P,\lambda)}\in \mathcal{D}'(X)$,
essentially the delta function determined by $(P,\lambda)$; here,
$\lambda$ is implicitly viewed by pull-back as a half-density on
$P$. To express this, given $\gamma \in
\mathcal{D}^\infty_{(1/2)}(X)$ write $\lambda=S_\lambda\cdot
\mathrm{dens}_P^{(1/2)}$ and $\gamma =T_\gamma \cdot
\mathrm{dens}_X^{(1/2)}$ for unique $S_\lambda\in
\mathcal{C}^\infty(P)$ and $T_\gamma\in \mathcal{C}^\infty(X)$; then
$\left <\delta _{(P,\lambda)},\gamma\right>=\int
_PS_\lambda\,T_\gamma \cdot \mathrm{dens}_P$.
Now $\delta _{(P,\lambda)}$ is a Lagrangian distribution: Let
$\phi=(p,q,\theta):U\rightarrow \mathbb{R}^\mathrm{n}\times
\mathbb{R}^\mathrm{n}\times \mathbb{R}$ be local coordinates for
$X$, centered at some $x_0\in P$, and defined on an open
neighborhood $U\ni x_0$. Suppose that $P\cap
U=\{p=\mathcal{P}(q),\theta=\Theta(q)\}\subseteq U$, where
$(\mathcal{P},\Theta):V=:\phi(U)\rightarrow
\mathbb{R}^\mathrm{n}\times \mathbb{R}$ is $\mathcal{C}^\infty$.
Then $(q)$ restricts to a system of local coordinates for $P$,
defined on $P\cap U$ and centered at $x_0$; accordingly, we shall
write $\mathrm{dens}_P^{(1/2)}=D_P\cdot \sqrt{\left |dq\right |}$,
for a unique $\mathcal{C}^\infty$ positive function $D_P$ on $P\cap
U$. Then if $\gamma$ is supported in $U$, we have
$\left <\delta _{(P,\lambda)},\gamma\right>=\int
_{\mathbb{R}^\mathrm{n}}S_\lambda(q)\,T_\gamma
\big(\mathcal{P}(q),q,\Theta(q)\big)\cdot
D_P(q)^2\,\left|dq\right|.$
On $\mathcal{C}^\infty_0(V)$, therefore, $\delta _{(P,\lambda)}$ is
the Fourier integral distribution
\begin{equation}\label{eqn:delta-P-lambda-local-fourier-integral}
\frac{1}{(2\pi)^{\mathrm{n}+1}}\int _{\mathbb{R}}\int
_{\mathbb{R}^{\mathrm{n}}}\,e^{i\ (\tau \,F+\eta \cdot H)}\,
S_\lambda(q)\, D_P(q)^2\,d\tau\,d\eta,
\end{equation}
where $F(p,q,\theta)=\theta-\Theta(q)$,
$H(p,q,\theta)=p-\mathcal{P}(q)$.

By its microlocal structure, the Szeg\"{o} projector of $X$ extends
to $\Pi:\mathcal{D}'(X)\rightarrow \mathcal{H}(X)$, where
$\mathcal{H}(X)\subseteq \mathcal{D}'(X)$ is the subspace of those
distributions all of whose Fourier components belong to the Hardy
space. Define
$\Delta:\mathbf{W}_{\mathrm{h}}\mathbf{P}_\mathfrak{I}\rightarrow
\mathcal{D}'(X)$ by $(P,\lambda)\mapsto\delta _{(P,\lambda)}$. Set
$u_{(P,\lambda)}=:\Pi \circ \Delta
(P,\lambda)=\Pi(\delta_{(P,\lambda)})\in \mathcal{H}(X)$, and for
 $k\in \mathbb{N}$ consider the
Fourier components, $u_{(P,\lambda),k}=:\Pi_k\circ \Delta
(P,\lambda)\in H(X)_k$; here $H(X)_k\cong H^0(M,A^{\otimes k})$ is
the level-$k$ Hardy space of $X$, and
$\Pi_k:\mathcal{D}'(X)\rightarrow H(X)_k$ the (extension of the)
$L^2$-orthogonal projector. This is the level-$k$ BPU map,
$\widetilde{\Phi} _k=:\Pi_k\circ
\Delta:\mathbf{W}_{\mathrm{h}}\mathbf{P}_\mathfrak{I}\rightarrow
H(X)_k$. Since by construction $\delta _{(P,\lambda)}$ is
$\mathbb{Z}_r$-invariant, so is $u_{(P,\lambda)}$; therefore
$u_{(P,\lambda),k}=0$, hence $\widetilde{\Phi} _k=0$, unless $r|k$.
Next suppose $k=l\cdot r$, $l\in \mathbb{N}$. By i) of Corollary 1.1
of \cite{dp}, $u_{(P,\lambda),k}(x)=O(k^{-\infty})$ whenever
$x\not\in S^1\cdot P$. If on the other hand $x\in S^1\cdot P$, by
ii) of the same Corollary in local Heisenberg coordinates for $X$
adapted to $P$ at $x$ and $\forall \,w\in T_mM\cong
\mathbb{C}^{\mathrm{n}}$, $m=:\pi(x)\in L\subseteq M$, there is an
asymptotic expansion:
\begin{eqnarray}
u_{(P,\lambda),k}\big(x+w/\sqrt k\big )&\sim&k^{\mathrm{n}/2}\cdot
r\cdot \left (\frac 2\pi\right )^{\mathrm{n}/2}\,e^{-ik\vartheta
(x)}\,S_\lambda (m)\,e^{-\|w^\perp\|^2-i\omega _m
(w^\perp,w^{\|})}\nonumber\\
&&+\sum _{f\ge 1}k^{(\mathrm{n}-f)/2}\,c_f(x,w).
\label{eqn:asymptotic-expansion-for-u}
\end{eqnarray}
Here, $w^{\|}\in T_mL$, $w^\perp\in (T_mL)^{\perp}$ denote the
orthogonal components of $w$, and $\omega _m (w^\perp,w^{\|})$ their
symplectic pairing. Furthermore, $e^{i\vartheta (x)}\in S^1$ is such
that $e^{i\vartheta (x)}\cdot x\in P$, and $\lambda=S_\lambda \cdot
\mathrm{dens}_L^{(1/2)}$, where $\mathrm{dens}_L^{(1/2)}$.

\begin{rem}
\label{rem:rapid-decay} By the arguments surrounding equations
(54)-(58), Lemma 3.7 and Claim 3.2 of \cite{dp},
$c_f(x,w)=\widehat{z}_f(x,w)\,e^{-\|w^\perp\|^2/2}$, where
$\widehat{z}_f(x,w)$ is a rapidly decaying function of $w^\perp$.
More precisely, up to some constant factor and oscillating term,
$\widehat{Z}_f(x,w)$ is the evaluation in $w^{\perp}$ of the Fourier
transform of the rapidly decreasing function $Z_j$ in statement ii)
of Lemma 3.7 of \cite{dp} (there $w^\perp=p_w$, $w^{\|}=q_w$).
\end{rem}

In particular, $\widetilde{\Phi} _k\big (P,\lambda\big)\neq 0$ if
$(P,\lambda)\in \mathbf{W}_{\mathrm{h}}\mathbf{P}_\mathfrak{I}$ and
$r|k$, $k\gg 0$. For $k=lr$, $l\in \mathbb{N}$, let
$\widetilde{\mathfrak{U}}_k\subseteq
\mathbf{W}_{\mathrm{h}}\mathbf{P}_\mathfrak{I}$ be the
$S^1$-invariant open subset where $\widetilde{\Phi} _k\neq 0$. Thus,
$\mathbf{W}_{\mathrm{h}}\mathbf{P}_\mathfrak{I}=\bigcup
_k\widetilde{\mathfrak{U}}_k$. Similarly, let
$\mathfrak{U}_k=:\widehat{\pi}(\widetilde{\mathfrak{U}}_k)=\widetilde{\mathfrak{U}}_k/S^1\subseteq
\mathbf{W}_{\mathrm{h}}\mathfrak{I}$; thus $\mathfrak{U}_k$ is open
in $\mathbf{W}_{\mathrm{h}}\mathfrak{I}$, and
$\mathbf{W}_{\mathrm{h}}\mathfrak{I}=\bigcup _{r|k}\mathfrak{U}_k$.

\begin{defn}\label{defn:proj-bpu-map}
For $k=lr$, $l\in \mathbb{N}$, define
$\Phi_k:\mathfrak{U}_k\rightarrow \mathbb{P}H(X)_k$ by
$(L,\lambda)\mapsto \left[\widetilde{\Phi} _k(P,\lambda)\right]$,
for any $P\in \mathbf{P}$ covering $L$. $\Phi_k$ is the
\textit{level-$k$ projectivized BPU map}.
\end{defn}


\section{The asymptotics of the differential of BPU maps}

Now we shall give an asymptotic expansion for certain scaling limits
of
$d_{(P,\lambda)}\widetilde{\Phi}_k:T_{(P,\lambda)}
\mathbf{W}_{\mathrm{h}}\mathbf{P}_\mathfrak{I}\rightarrow
H(X)_k,$
at a given $(P,\lambda)\in
\mathbf{W}_{\mathrm{h}}\mathbf{P}_\mathfrak{I}$. We may assume
without loss that $k=r\cdot l$, $l\in \mathbb{N}$. Assuming that
$\Delta$ is differentiable, since $\Pi _k$ is linear we have
\begin{equation}\label{eqn:diff-linear-non-linear}
 d_{(P,\lambda)}\widetilde{\Phi}_k\big(\widetilde{W}\big)=\Pi _k\Big
 (d_{(P,\lambda)}\Delta\big(\widetilde{W}\big)\Big)\,\,\,\,\,\,\,\,\,\,\,\,(\widetilde{W}\in T_{(P,\lambda)}
\mathbf{W}_{\mathrm{h}}\mathbf{P}_\mathfrak{I}).
\end{equation}
We shall first determine
$d_{(P,\lambda)}\Delta\big(\widetilde{W}\big)$. To this end, set
$L=\pi(P)$, and suppose $\widetilde{W}=\widetilde{W}(f,\ell)$, with
$\int _L\ell \bullet \lambda=0$ (Proposition
\ref{prop:tg-space-half-weights}). Identify $\lambda$ and $\ell$
with their pull-backs to $P$, and write $\lambda=S_\lambda\cdot
\mathrm{dens}_P^{(1/2)}$ and $\ell=S_\ell\cdot
\mathrm{dens}_P^{(1/2)}$; the smooth functions $S_\lambda,\,S_\ell$
on $P$ descend on $L$. Locally near some $x_0\in P$, fix good local
coordinates $(p,q,\theta)$ for $X$ along $P$ centered at $x_0$,
defined on $X'\subseteq X$ (Definition
\ref{defn:good-local-coordinates}). Then $(p,q)$ are naturally local
coordinates for $M$ centered at $m_0=:\pi(x_0)$, defined on
$M'=:\pi(X')$; the projection $X'\rightarrow M'$ is represented by
$(p,q,\theta)\mapsto (p,q)$. Set $P'=:X'\cap P=\{p=0,\,\theta =0\}$,
$L'=L\cap M'=:\{p=0\}$. We may view $(q)$ as local coordinates on
$P'$. Perhaps after restricting $X'$,
$\left.\pi\right|_{P'}:P'\rightarrow L'$ is an isometric
diffeomorphism.

\begin{prop}
\label{prop:derivative-of-Delta}
In the notation of the preceding discussion, the following holds:

\begin{enumerate}
  \item Let us extend $f$ to some tubular neighborhood of $L$ by the normal cotangent structure,
  and let $\upsilon_f$ be its Hamiltonian vector field. Then, $\forall\,m\in L'$, we have
$\upsilon_f(m)= \sum _{j=1}^\mathrm{n}a_j(m)\cdot
\left.\frac{\partial}{\partial p_j}\right |_m$,
for unique $a_j\in \mathcal{C}^\infty(L')$.
  \item  Let $\Gamma (L,f)$ be as in Definition
\ref{derivative-of-density-along-flow}, $a=(a_j)$. Locally near
$x_0$, $d_{(P,\lambda)}\Delta\big(\widetilde{W}\big)\in
\mathcal{D}'(X)$ is the Fourier integral
$$
\frac{1}{(2\pi)^{\mathrm{n}+1}}\,\int
_{\mathbb{R}}\!\int_{\mathbb{R}^\mathrm{n}}\!e^{i\, (\tau
\theta+\eta\cdot p)}\,\Big [\Big (S_\ell+S_\lambda\cdot\Gamma
(L,f)\Big)-i\Big(\tau\,f+\eta\cdot
a\Big)\,S_\lambda\Big]\,D_{P}^2\,d\tau\,d\eta,
$$
where $\mathrm{dens}_P^{(1/2)}=D_P\cdot \sqrt{|dq|}$ is the
Riemannian half-density on $P$ (or $L$).
\end{enumerate}
\end{prop}

\textit{Proof of Proposition \ref{prop:derivative-of-Delta}.} By the
discussion surrounding Proposition \ref{prop:tg-space-weights},
$\upsilon_f(m)\in (T_mL)^\perp=J_m(T_mL)$, $\forall\,m\in L$. This
proves 1., in view of Proposition
\ref{prop:local-adapted-coordinates}.

For some $\epsilon>0$, suppose $\gamma
:(-\epsilon,\epsilon)\rightarrow
\mathbf{W}_{\mathrm{h}}\mathbf{P}_\mathfrak{I}$, $\gamma
(t)=(P_t,\lambda_t)$, is $\mathcal{C}^\infty$ with $\gamma
(0)=(P,\lambda)$, $\gamma'(0)=\widetilde{W}$. Then $L_t=:\pi(P_t)\in
\mathfrak{I}$ $\forall\,t\in (-\epsilon,\epsilon)$. If $\phi _t$ is
the local flow of $\upsilon_f$, then $L_t=\phi _t(L)$ to first order
in $t$. Let $\upsilon_f^\sharp$ be horizontal lift of $\upsilon_f$
to $X$, and $\frac{\partial}{\partial \theta}$ the generator of the
$S^1$-action. Then $\widetilde{\upsilon}_f=:\upsilon_f^\sharp-f\cdot
\frac{\partial}{\partial \theta}$ is a contact vector field on
$\pi^{-1}(M')\supseteq P$, whose local flow $\widetilde{\phi} _t$
covers $\phi_t$, and $P_t=\widetilde{\phi} _t(P)$ to first order in
$t$. Next, let $\beta_t:L_t\rightarrow L$ be induced by the normal
cotangent structure near $L$. There is a smooth path $\eta _t$ of
half-weights on $L$, such that $\lambda _t=\beta_t^*(\eta_t)$ for
every $t$. Then $\ell=\eta'(0)$, so that $\eta _t=\lambda +t\cdot
\ell +O(t^2)$. By 1., we have $\widetilde{\upsilon}_f(x)= \sum
_{j=1}^\mathrm{n}a_j(q)\cdot \left.\frac{\partial}{\partial
p_j}\right |_x-f(q)\cdot \left.\frac{\partial}{\partial
\theta}\right|_x$, if $x\in P'$ has local coordinates $(0,q,0)$.
Thus, for $t\sim 0$, $P_t\subseteq X$ is locally defined by
$p_j=t\cdot a_j(q)+O(t^2)$ ($j=1,\ldots,\mathrm{n}$), and $\theta
=-t\cdot f(q)+O(t^2)$. Write $\lambda _t=S_{\lambda_t}\cdot
\mathrm{dens}_{L_t}^{(1/2)}$, for unique $S_{\lambda _t}\in
\mathcal{C}^\infty(L_t)$. By
(\ref{eqn:delta-P-lambda-local-fourier-integral}), $\delta
_{(P_t,\lambda_t)}=\Delta \big (P_t,\lambda _t\big)$ is locally near
$x_0$ the Fourier integral
\begin{equation}
\label{eqn:fourier-integral-for-delta}
\frac{1}{(2\pi)^{\mathrm{n}+1}}\,\int
_{\mathbb{R}}\!\int_{\mathbb{R}^\mathrm{n}}\!e^{i\, (\tau
F_t+\eta\cdot H_t)}\,S_{\lambda_t}(q)\,D_{P_t}(q)^2\,d\tau\,d\eta,
\end{equation}
where $F_t(p,q,\theta)=\theta-tf(q)+O(t^2)$,
$H_t(p,q,\theta)=p-t\,a(q)+O(t^2)$. Here $p=(p_j)$, $a=(a_j)$ (cfr
Lemma 2.2 of \cite{dp}). Furthermore, $q=(q_j)$ restrict to local
coordinates on $L_t$, and $S_{\lambda_t}(q)$, $D_{P_t}(q)$ are meant
in this local coordinate system (thus,
$\mathrm{dens}_{P_t}^{(1/2)}=D_{P_t}\cdot \sqrt{|dq|}$). By \S 7.8
of \cite{hor}, the $t$-derivative of $\Delta\big (P_t,\lambda_t\big
)$ may be computed by differentiating with respect to $t$ under the
integral sign in (\ref{eqn:fourier-integral-for-delta}).

\begin{lem}\label{lem:local-representation-for-beta_t}
\textbf{A}: When $q$ is adopted as a system of local coordinates for
both $L=L_0$ and $L_t$, we have $\beta _t(q)=q+O(t^2)$. \textbf{B}:
Let $\left .\phi_t\right |_L:L\rightarrow L_t=\phi _t(L)$ be the
diffemomorphism induced by the flow of $\upsilon _f$. Let
$\left(\left .\phi_t\right |_L\right )^{-1}:L_t\rightarrow L$ be the
inverse diffemorphism. Then $\left(\left .\phi_t\right |_L\right
)^{-1}=\beta _t+O(t^2)$.
\end{lem}

\textit{Proof.} Let $U\subseteq M$ be the open tubular neighborhood
of $L$ produced in the construction of the normal cotangent
structure, and let $\widehat{\pi}:U\rightarrow L$ be the normal
cotangent projection. Let $\widehat{\pi}':U\rightarrow L$ be the
locally defined projection $M'\rightarrow L'$ which is given in good
local coordinates by $(p,q)\mapsto q$. By the properties of good
local coordinates, the fibers of both $\widehat{\pi}$ and
$\widehat{\pi}'$ meet $L$ perpendicularly at each $m\in L'$. It
follows that, in local coordinates,
$\widehat{\pi}(m)-\widehat{\pi}'(m)=O\left(\mathrm{dist}(m,L)^2\right)$
($m\in M'$). Restricting to $P_t$, this implies \textbf{A}, and the
proof of \textbf{B} is similar. \hfill $\square$


\medskip

As in the previous discussion, let $\eta_t$ be the half-weight on
$L$ such that $\beta _t^*(\eta_t)=\lambda_t$. Let us write
$\eta_t=S_{\eta_t}\cdot \mathrm{dens}_L^{(1/2)}$, $\ell =S_\ell
\cdot \mathrm{dens}_L^{(1/2)}$ for uniquely determined
$\mathcal{C}^\infty$ functions $S_{\eta_t}$ and $S_\ell$ on $L$.
Therefore, $S_{\eta_t}=S_\lambda+t\,S_\ell+O(t^2)$. Notice that
$S_{\lambda_t}$ is a smooth function on $P_t$, while $S_{\eta_t}$ is
a smooth function on $P=P_0$. Since $q$ restricts to a system of
local coordinates on both $P$ and $P_t$, $t\sim 0$, we can consider
the local expressions $S_{\lambda _t}(q)$ and $S_{\eta_t}(q)$. By
Definition \ref{derivative-of-density-along-flow} and Lemma
\ref{lem:local-representation-for-beta_t},
$D_{P_t}(q)/D_P(q)=1+t\,\Gamma (L,f)+O(t^2)$.
In view of Corollary \ref{lem:local-representation-for-beta_t},
$\beta _t^*(\sqrt{|dq|})=\sqrt{|dq|}+O(t^2)$, and
$(\beta_t^*g)(q)=g(q)+O(t^2)$ for every locally defined function
$g$. Therefore,
\begin{eqnarray*}
\lefteqn{\lambda_t=\beta_t^*\big(S_{\eta_t}\cdot
\mathrm{dens}_P\big)=\beta_t^*\Big(S_{\eta_t}\,D_P\cdot
\sqrt{|dq|}\Big)= S_{\eta_t}\,D_P\cdot\sqrt{|dq|}+O(t^2)} \\
&&= S_{\eta_t}\,\frac{D_P}{D_{P_t}}\,D_{P_t}\cdot\sqrt{|dq|}+O(t^2)=
S_{\eta_t}\,\frac{D_P}{D_{P_t}}\cdot \mathrm{dens}_{P_t}+O(t^2).
\end{eqnarray*}
We deduce that
\begin{eqnarray}
\label{eqn:S-lambda-t}
S_{\lambda_t}&=&S_{\eta_t}\,\frac{D_P}{D_{P_t}}+O(t^2)=\Big(S_\lambda+tS_\ell\Big
)\cdot \Big (1-t\,\Gamma (L,f)\Big )+O(t^2)\nonumber\\
&=&S_\lambda+t\Big (S_\ell-S_\lambda\cdot\Gamma
(L,f)\Big)+O(t^2),\nonumber
\end{eqnarray}
whence
$S_{\lambda_t}\cdot D_{P_t}^2=D_P^2\cdot \Big[S_\lambda+t\Big
(S_\ell+S_\lambda\cdot \Gamma (L,f)\Big)\Big]+O(t^2).
$
The proof of the second statement of Proposition
\ref{prop:derivative-of-Delta} is completed by inserting this
equality
in
(\ref{eqn:fourier-integral-for-delta}), and differentiating with
respect to $t$ under the integral sign at $t=0$.
\hfill $\square$

\bigskip

Thus, locally near $x_0$, we have
$d_{(P,\lambda)}\Delta \big(\widetilde{W}\big) =\sum
_{j=1}^4d_{(P,\lambda)}\Delta \big(\widetilde{W}\big)_j$,
where
\begin{equation}\label{eqn:first-piece}
d_{(P,\lambda)}\Delta
\big(\widetilde{W}\big)_j=\frac{1}{(2\pi)^{\mathrm{n}+1}}\,\int
_{\mathbb{R}}\!\int_{\mathbb{R}^\mathrm{n}}\!e^{i\, (\tau
\theta+\eta\cdot p)}\,b_j\,D_P^2\,d\tau\,d\eta,
\end{equation}
with $b_1=:S_\ell$, $b_2=:S_\lambda\cdot \Gamma(L,f)$, $b_3=:-i\tau
f S_\lambda$, $b_4=:-i(\eta \cdot a)\, S_\lambda$. Applying the
level-$k$ Szeg\"{o} kernel, by (\ref{eqn:diff-linear-non-linear}) we
obtain
$ d_{(P,\lambda)}\widetilde{\Phi}_k\big(\widetilde{W}\big)=\sum
_{j=1}^4
d_{(P,\lambda)}\widetilde{\Phi}_k\big(\widetilde{W}\big)_j$,
where
$d_{(P,\lambda)}\widetilde{\Phi}_k\big(\widetilde{W}\big)_j=:\Pi_k\Big
(d_{(P,\lambda)}\Delta \big(\widetilde{W}\big)_j\Big)$, $1\le j\le
4$. Let us now consider the transverse scaling asymptotics for
$d_{(P,\lambda)}\widetilde{\Phi}_k\big(\widetilde{W}\big)_j$ near
$x_0$. Given $x\in P'$ with good local coordinates $(0,q,0)$, and
$w\in \mathbb{R}^\mathrm{n}$, $x+w$ will mean the point in $P'$
having good local coordinates $(w,q,0)$. The real n-space
$\mathbb{R}^\mathrm{n}\subseteq \mathbb{C}^\mathrm{n}$ is unitarily
identified with the orthocomplement $\left (T_mL\right )^\perp$,
$m=:\pi (x)$, hence with a subspace of $\mathrm{Hor}(X)_x\subseteq
T_xX$.

\begin{lem}\label{lem:scaling-asympt-first-piece}
Suppose $x_0\in P'\subseteq P$ is a sufficiently small open
neighborhood. Uniformly in $x\in S^1\cdot P'$ and $w \in
\mathbb{R}^{\mathrm{n}}$ of bounded norm, for $j=1,2$ the following
asymptotic expansion holds as $l\rightarrow +\infty$ and $k=l\cdot
r$:
\begin{eqnarray}\label{eqn:first-asymptotic-expansion}
d_{(P,\lambda)}\widetilde{\Phi}_k\big(\widetilde{W}\big)_j\left(x+\frac{w}{\sqrt{k}}\right
)&\sim&k^{\mathrm{n}/2}\cdot r\, \left (\frac 2\pi\right
)^{\mathrm{n}/2}\, e^{-ik\vartheta(x)}\,b_j (m)\, e^{-\|w\|^2}\nonumber \\
&&+\sum _{h\ge 1}C_{hj}(x,w)\, k^{(\mathrm{n}-h)/2},
\end{eqnarray}
where $m=:\pi(x)\in L$, and $\vartheta(x)\in (-\pi,\pi]$ is such
that $e^{i\vartheta}\cdot x\in P'$. Similarly,
\begin{eqnarray}\label{eqn:second-asymptotic-expansion}
d_{(P,\lambda)}\widetilde{\Phi}_k\big(\widetilde{W}\big)_3\left(x+\frac{w}{\sqrt{k}}\right
)&\sim&-i\, k^{1+\mathrm{n}/2}\cdot r\, \left (\frac 2\pi\right
)^{\mathrm{n}/2}\,  e^{-ik\vartheta(x)}\,f (m)\,S_\lambda(m)\,
e^{-\|w\|^2}\nonumber \\
&&+\sum _{h\ge 1}C_{h3}(x,w)\, k^{1+(\mathrm{n}-h)/2},
\end{eqnarray}
\begin{eqnarray}\label{eqn:third-asymptotic-expansion}
d_{(P,\lambda)}\widetilde{\Phi}_k\big(\widetilde{W}\big)_4\left(x+\frac{w}{\sqrt{k}}\right
)&\sim&\sum _{h\ge 0}C_{h4}(x,w)\,k^{(1+\mathrm{n}-h)/2}.
\end{eqnarray}
Furthermore, for $j=1,2,3,4$ and $h\ge 1$ ,
$C_{hj}(x,w)=\widehat{C}_h^{(j)}(x,w)\,e^{-\|w\|^2/2},$ where
$\widehat{C}_h^{(j)}$ is a rapidly decaying function of $w$.
\end{lem}

\textit{Proof.} We may assume $x\in P$. For $j=1,2$, we have
$d_{(P,\lambda)}\widetilde{\Phi}_k\big(\widetilde{W}\big)_j=\Pi
_k\big (\delta _{(P,\sigma _j)}\big)$, where $\sigma _j$ is the
half-density (not necessarily a half-weight) on $L$ defined by
$\sigma _j=b_j\cdot \mathrm{dens}_L^{(1/2)}$. By Corollary 1.1 of
\cite{dp}, the scaling asymptotics of $\Pi _k\big (\delta
_{(P,\sigma _j)}\big)$ at $x\in P$ are given - in Heisenberg local
coordinates adapted to $P$ at $x$ - by asymptotic expansions akin to
(\ref{eqn:asymptotic-expansion-for-u}), with $b_j$ in place of
$S_\lambda$. In (\ref{eqn:asymptotic-expansion-for-u}), $w$ is
allowed to vary in $\mathbb{C}^\mathrm{n}$, and $x+w$ denotes the
point of $X$ having adapted Heisenberg local coordinates $(w,0)$. On
the other hand, good local coordinates along $P$ are constructed by
glueing moving systems of trasverse Heisenberg local coordinates
along a system of arbitrary local coordinates along $P$ (this is
made precise in Proposition \ref{prop:local-adapted-coordinates}).
Since in (\ref{eqn:first-asymptotic-expansion}) $w$ is required to
be a real vector, the expression $x+\frac{w}{\sqrt{k}}$ (in the
given system of good local coordinates) represents a transverse
displacement from $P$ which is also represented by the expression
$x+\frac{w}{\sqrt{k}}$ in a system of Heisenberg local coordinates
adapted to $P$ at $x$. Thus, an expansion of type
(\ref{eqn:asymptotic-expansion-for-u}) still holds in good local
coordinates, so far as the rescaling occurs in the transverse
direction only.

The proof of (\ref{eqn:second-asymptotic-expansion}) is similar, but
we need to explain the extra factor of $k$. To this end, we remark
that $d_{(P,\lambda)}\widetilde{\Phi}_k\big(\widetilde{W}\big)_3=\Pi
_k\Big (d_{(P,\lambda)}\Delta\big(\widetilde{W}\big)_3\Big)$. Now
$d_{(P,\lambda)}\Delta(W)_3$ is the Fourier integral
(\ref{eqn:first-piece}), with $b_3=-i\tau fS_\lambda$ (introduce a
cut-off to make this compactly supported near $x$). Due to the
factor $\tau$ appearing in the amplitude, this is not of the form
$\delta_{(P,\sigma)}$ for a $\mathcal{C}^\infty$ half-density on
$P$. However, the techniques in the proof of
(\ref{eqn:asymptotic-expansion-for-u}) can still be applied. Namely,
one applies to (\ref{eqn:first-piece}) the Boutet de Monvel -
Sj\"{o}strand parametrix for the Szeg\"{o} kernel, and then takes
the $k$-th Fourier component. After suitably rescaling the
integration variables involved, this yields an oscillatory integral
to which the stationary phase Lemma may be applied. By Claim 3.2 of
\cite{dp}, this leads to a unique stationary point where $\tau=1$
and $\eta=0$. The rescaling involved in $\tau$ is $\tau\mapsto
k\,\tau$, and this accounts for the extra factor of $k$ in
(\ref{eqn:second-asymptotic-expansion}).

Let us consider (\ref{eqn:third-asymptotic-expansion}). Now
$d_{(P,\lambda)}\widetilde{\Phi}_k\big(\widetilde{W}\big)_4=\Pi
_k\Big (d_{(P,\lambda)}\Delta\big(\widetilde{W}\big)_4\Big)$, and
$d_{(P,\lambda)}\Delta\big(\widetilde{W}\big)_4$ is the Fourier
integral (\ref{eqn:first-piece}), with $b_4=-i(\eta \cdot
a)S_\lambda$. The same arguments used in the previous paragraph
apply, so that $\Pi _k\Big
(d_{(P,\lambda)}\Delta\big(\widetilde{W}\big)_4\Big)\left(x+\frac{w}{\sqrt{k}}\right
)$ is an oscillatory integral to which the stationary phase Lemma
may be applied. By the arguments in the proof of Theorem 1.1 of
\cite{dp}, the rescaling in $\eta$ is $\eta\mapsto k^{3/2}\eta$,
hence the leading order term of the resulting asymptotic expansion
has degree at most $k^{(\mathrm{n}+3)/2}$. However, as mentioned the
stationary point of the phase has $\eta=0$. Since by Theorem 7.7.5
of \cite{hor} the first term involving derivatives of the amplitude
has degree $(3+\mathrm{n})/2-1=(1+\mathrm{n})/2$, the terms in
$k^{(\mathrm{n}+3)/2}$ and $k^{1+\mathrm{n}/2}$ both vanish.

The last statement is proved arguing as in Remark
\ref{rem:rapid-decay}. \hfill $\square$

\bigskip

The $\mathcal{C}^\infty$ $\mathbb{R}^\mathrm{n}$-valued function $a$
on $P$ appearing in $b_4$ depends linearly on $f$ and is independent
of $\ell$, while $S_\ell$ depends linearly on $\ell$, and is
independent of $f$. Let us write $W=W(f,\ell)=W(f,0)+W(0,\ell)\in
T_{(P,\lambda)}\mathbf{W}_{\mathrm{h}}\mathbf{P}_\mathfrak{I}$. With
$x\in S^1\cdot P$, $m=\pi(x)$ and $w\in \mathbb{R}^\mathrm{n}$, we
have
\begin{eqnarray}
\label{eqn:cotribution-f-0} \lefteqn{
d_{(P,\lambda)}\widetilde{\Phi}_k\Big (\widetilde{W}(f,0)\Big)\left
(x+\frac{w}{\sqrt{k}}\right )
=\sum _{j=2}^4d_{(P,\lambda)}\widetilde{\Phi}_k
\big(\widetilde{W}\big)_j}\\
&&\sim -i\, k^{1+\mathrm{n}/2}\cdot r\, \left (\frac 2\pi\right
)^{\mathrm{n}/2}\,  e^{-ik\vartheta(x)}\,f (m)\,S_\lambda(m)\,
e^{-\|w\|^2}
+\sum _{h\ge 1}D_{h}(x,w)\, k^{1+(\mathrm{n}-h)/2}\nonumber,
\end{eqnarray}
\begin{eqnarray}
\label{eqn:cotribution-l-0} \lefteqn{
d_{(P,\lambda)}\widetilde{\Phi}_k\Big
(\widetilde{W}(0,\ell)\Big)\left (x+\frac{w}{\sqrt{k}}\right
)=d_{(P,\lambda)}\widetilde{\Phi}_k \big(\widetilde{W}\big)_2 }
\\
&&\sim  k^{\mathrm{n}/2}\cdot r\, \left (\frac 2\pi\right
)^{\mathrm{n}/2}\, e^{-ik\vartheta(x)}\,S_\ell (m)\, e^{-\|w\|^2}
+\sum _{h\ge 1}E_{h}(x,w)\, k^{(\mathrm{n}-h)/2}.\nonumber
\end{eqnarray}

For $\widetilde{W}=\widetilde{W}(f,\ell)\in
T_{(P,\lambda)}\mathbf{W}_{\mathrm{h}}\mathbf{P}_\mathfrak{I}$ and
$k=1,2,\ldots$, set $\widetilde{W}_k=:\widetilde{W}(f,k\,\ell)$.
Thus, $d_{(P,\lambda)}\widetilde{\Phi}_k\big
(\widetilde{W_k}\big)=d_{(P,\lambda)}\widetilde{\Phi}_k\big
(\widetilde{W}(f,0)\big)+k\,d_{(P,\lambda)}\widetilde{\Phi}_k\big
(\widetilde{W}(0,\ell)\big)$. Given Lemma
\ref{lem:scaling-asympt-first-piece}, summing over $j$ we obtain:

\begin{cor}\label{cor:rescaled-asympt-w-k}
Suppose $(P,\lambda)\in
\mathbf{W}_{\mathrm{h}}\mathbf{P}_\mathfrak{I}$, $\widetilde{W}\in
T_{(P,\lambda)}\mathbf{W}_{\mathrm{h}}\mathbf{P}_\mathfrak{I}$.
Then:
\begin{itemize}
  \item If $x\not\in S^1\cdot P$, then
$d_{(P,\lambda)}\widetilde{\Phi}_k\big(\widetilde{W}_k\big)(x)=O(k^{-\infty})$,
uniformly in $x$ on compact subsets of $X\setminus S^1\cdot P$.
  \item Uniformly in $x\in S^1\cdot P$ and in $w\in T_{\pi(x)}L^\perp\subseteq
  T_{\pi(x)}M$ of bounded norm, the
following asymptotic expansion holds as $l\rightarrow +\infty$ and
$k=l\cdot r$:
\begin{eqnarray*}
d_{(P,\lambda)}\widetilde{\Phi}_k\big(\widetilde{W}_k\big)\left(x+\frac{w}{\sqrt{k}}\right
)&\sim& k^{1+\mathrm{n}/2}\cdot r\, \left (\frac 2\pi\right
)^{\mathrm{n}/2}\, e^{-ik\vartheta(x)}\, \gamma _{\ell f}(m)\,
e^{-\|w\|^2}\nonumber \\
&&+\sum _{h\ge 1}H_h(x,w)\, k^{1+(\mathrm{n}-h)/2},
\end{eqnarray*}
where $\gamma_{\ell f}=S_\ell-i\,f \,S_\lambda$, and $\forall\,h\ge
1$ we have $H_h(x,w)=\widehat{H}_h(x,w)\,e^{-\|w\|^2/2}$,
$\widehat{H}_h$ being a rapidly decaying function of $w$.
\end{itemize}
\end{cor}

We can now prove:
\begin{prop}\label{prop:asymptotically-orthogonal}
If $(P,\lambda)\in \mathbf{W}_{\mathrm{h}}\mathbf{P}_\mathfrak{I}$,
as $l\rightarrow +\infty$ and $k=l\cdot r$ we have:
\begin{eqnarray}\label{eqn:asympt-exp-norm-phi-k}
\left
(\widetilde{\Phi}_k(P,\lambda),\widetilde{\Phi}_k(P,\lambda)\right)_{L^2(X)}
&\sim&k^{\mathrm{n}/2}\,\left (\frac{2}{\pi}\right
)^{\mathrm{n}/2}\,r^2+\sum _{h\ge 1}q_h\cdot k^{(\mathrm{n}-h)/2}.
\end{eqnarray}
Given tangent vectors
$\widetilde{W}=\widetilde{W}(f,\ell),\,\widetilde{W}'=\widetilde{W}(f',\ell')\in
T_{(P,\lambda)}\mathbf{W}_{\mathrm{h}}\mathbf{P}_\mathfrak{I}$, set
$$F\big(\widetilde{W},\widetilde{W}'\big)=:\big
(S_\ell\,S_{\ell'}+f\,f'\,S_\lambda^2\big)+i\big(S_\ell\,f'-f\,S_{\ell'}\big)\,S_\lambda.
$$
Then the following asymptotic expansions hold as $l\rightarrow
+\infty$ and $k=l\cdot r$:
\begin{eqnarray}\label{eqn:asympt-exp-hermitian-product-phi-k}
\left (d_{(P,\lambda)}\widetilde{\Phi}_k\big(\widetilde{W}_k\big),
d_{(P,\lambda)}\widetilde{\Phi}_k\big(\widetilde{W}'_k\big)\right)_{L^2(X)}&\sim&k^{2+\mathrm{n}/2}\cdot
r^2\left (\frac{2}{\pi}\right
 )^{\mathrm{n}/2}\int _LF\big(\widetilde{W},\widetilde{W}'\big)\cdot
 \mathrm{dens}_L\nonumber\\
&&+\sum _{h\ge 1}r_{h}k^{2+(\mathrm{n}-h)/2},
\end{eqnarray}
\begin{eqnarray}\label{eqn:asympt-exp-hermitian-product-phi-k-dphi-k}
\left
(\widetilde{\Phi}_k(P,\lambda),d_{(P,\lambda)}\widetilde{\Phi}_k\big(\widetilde{W}_k\big)\right)_{L^2(X)}&\sim&
k^{1+\mathrm{n}/2}\cdot i\,r^2\,\left (\frac 2\pi\right
)^{\mathrm{n}/2}\,\int _LS_\lambda ^2\,f\cdot
\mathrm{dens}_L\nonumber \\
&&+\sum _{h\ge 1}s_hk^{1+(\mathrm{n}-h)/2}.
\end{eqnarray}
\end{prop}
An estimate similar to (\ref{eqn:asympt-exp-norm-phi-k}) was first
proved in \cite{bpu}, where BPU maps where originally phrased using
Fourier-Hermite distributions and symplectic spinors.

\begin{rem}
\label{rem:integral-of-F-is-product} Suppose $(P,\lambda)\in
\mathbf{P}_\mathfrak{I}$, and set $L=\pi(P)$. If
$W=W(f,\ell),\,W'=W'(f',\ell')\in
T_{(L,\lambda)}\mathbf{W}_{\mathrm{h}}\mathfrak{I}$, and
$\widetilde{W}=\widetilde{W}(f,\ell),\,\widetilde{W}'=\widetilde{W}(f',\ell')
\in  T_{(P,\lambda)}\mathbf{W}_{\mathrm{h}}\mathbf{P}_\mathfrak{I}$
are their lifts, then $ \int
_L\,F\big(\widetilde{W},\widetilde{W}'\big)\cdot
\mathrm{dens}_L=G_{(L,\lambda)}\big(W,W'\big)
+i\,\Omega_{(L,\lambda)}\big(W,W'\big).$
\end{rem}

\textit{Proof of Proposition \ref{prop:asymptotically-orthogonal}.}
Since $\Phi _k(P,\lambda)=\Pi _k(\delta_{P,\lambda})$, arguing as in
(75) of \cite{dp} we have
$
\left (\Phi _k(P,\lambda),\Phi _k(P,\lambda)\right )_{L^2(X)}=\left
<\delta_{P,\lambda},\overline{\Phi _k(P,\lambda)}\right
>=\int _PS_\lambda\, \overline{\Phi _k(P,\lambda)}\cdot
\mathrm{dens}_P$.
Now (\ref{eqn:asymptotic-expansion-for-u}) with $w=0$ yields an
asymptotic expansion for $\Phi _k(P,\lambda)(x)$, $x\in P$; when
$x\in P$, we may actually assume $\vartheta(x)=0$ in
(\ref{eqn:asymptotic-expansion-for-u}). Inserting the latter
asymptotic expansion in
the former integral proves (\ref{eqn:asympt-exp-norm-phi-k}), since
$\int _PS_\lambda ^2\cdot \mathrm{dens}_P=r \int _L S_\lambda^2\cdot
\mathrm{dens}_L =r\int _L\lambda\bullet \lambda=r,$ because
$P\rightarrow L$ is a Riemannian covering of degree $r$, and
$\lambda=S_\lambda \cdot \mathrm{dens}_L^{(1/2)}$ is a half-weight
on $L$. The proof of
(\ref{eqn:asympt-exp-hermitian-product-phi-k-dphi-k}) is similar,
except that we now need to use the asymptotic expansion in Corollary
\ref{cor:rescaled-asympt-w-k}, and recall that $\int
_LS_\lambda\,S_\ell\cdot \mathrm{dens}_L=\int _L\lambda \bullet \ell
=0$.

Let us now consider (\ref{eqn:asympt-exp-hermitian-product-phi-k}).
Let $U\supseteq L=\pi(P)$ be a suitably small tubular neighborhood
of $L$ in $M$, so that $T=:\pi^{-1}(U)\subseteq X$ is an
$S^1$-invariant open neighborhood of $P$. In view of the first
statement of Corollary \ref{cor:rescaled-asympt-w-k}, we have:
\begin{equation}\label{eqn:hermitian-product-asymptotic-derivative-tubular}
\left (d_{(P,\lambda)}\widetilde{\Phi}_k\big(\widetilde{W}_k\big),
d_{(P,\lambda)}\widetilde{\Phi}_k\big(\widetilde{W}'_k\big)\right)_{L^2(X)}\sim\int
_Td_{(P,\lambda)}\widetilde{\Phi}_k\big(\widetilde{W}_k\big)\,
\overline{d_{(P,\lambda)}\widetilde{\Phi}_k\big(\widetilde{W}'_k\big)}\cdot
\mathrm{dens}_X.
\end{equation}

Suppose that $U=\bigcup _jU_j$ is an open cover of $U$, such that on
each $T_j=:\pi^{-1}(U_j)$ there is a system of good local
coordinates for $X$ near $P$, in the stronger sense of Proposition
\ref{prop:local-adapted-coordinates-planckian}; we may as well
assume that the $T_j$'s are finitely many. Let $\{\varphi _j\}$ be a
partition of unity on $U$ subordinate to the open cover $\{U_j\}$,
so that $\{\widetilde{\varphi} _j\}$ is a partition of unity on $T$
for the open cover $\{T_j\}$, where $\widetilde{\varphi} _j=:\varphi
_j\circ \pi$. Given
(\ref{eqn:hermitian-product-asymptotic-derivative-tubular}), we
obtain
$\left (d_{(P,\lambda)}\widetilde{\Phi}_k\big(\widetilde{W}_k\big),
d_{(P,\lambda)}\widetilde{\Phi}_k\big(\widetilde{W}'_k\big)\right)_{L^2(X)}\sim
\sum_jA_{jk}$,
where we have set $A_{jk}=:\int
_T\widetilde{\varphi}_j\,d_{(P,\lambda)}\widetilde{\Phi}_k\big(\widetilde{W}_k\big)\,
\overline{d_{(P,\lambda)}\widetilde{\Phi}_k\big(\widetilde{W}'_k\big)}\cdot
\mathrm{dens}_X$. Let us now evaluate each $A_{jk}$ asymptotically
as $k\rightarrow +\infty$.

Let us set
$R_k=:d_{(P,\lambda)}\widetilde{\Phi}_k\big(\widetilde{W}_k\big)\,
\overline{d_{(P,\lambda)}
\widetilde{\Phi}_k\big(\widetilde{W}'_k\big)}$. In good local
coordinates,
$A_{jk}=\int _{T_j}\varphi _j\cdot R_k(p,q,\theta)\,
D_P\,dp\,dq\,d\theta=k^{-\mathrm{n}/2}\int \int \varphi _j \cdot
R_k\cdot D_X^2\left (\frac{w}{\sqrt
k},q,\theta\right)\,dw\,dq\,d\theta$,
where we have performed the rescaling $p=\frac{w}{\sqrt{k}}$, and
written $\mathrm{dens}_X=\frac{1}{2\pi}\,D_X^2\cdot
|dp\,dq\,d\theta|$. Clearly, $\left (\frac{w}{\sqrt
k},q,\theta\right)$ corresponds to $x+\frac{w}{\sqrt{k}}$, where
$x\in \big(S^1\cdot P\big)\cap T_j$ has good local coordinates
$(0,q,\vartheta)$. On the other hand, by the second statement of
Corollary \ref{cor:rescaled-asympt-w-k} working in good local
coodinates we have
\begin{equation}\label{eqn:asympt-exp-product-diff}
R_k\left(x+\frac{w}{\sqrt k}\right )\sim k^{2+\mathrm{n}}\cdot
r^2\,\left (\frac{2}{\pi}\right
)^{\mathrm{n}}F(W,W')\,e^{-2\|w\|^2}+\sum _{h\ge
1}t_h(x,w)k^{2+\mathrm{n}-h/2},
\end{equation}
where for every $h\ge 1$ we have
$t_h(x,w)=\widehat{t}_h(x,w)\,e^{-\|w\|^2}$, with $\widehat{t}_h$ a
rapidly decaying function of $w$. Now we can perform the integration
over $T_j$ by first integrating in $dw$ over $\mathbb{R}^\mathrm{n}$
and in $\frac{1}{2\pi}\,d\theta$ over $(-\pi,\pi)$, and then in $dq$
(viewed now as a system of local coordinates on $L$). To perform the
former, we remark that by the construction of good local
coordinates, we have
$D_X\left(x+\frac{w}{\sqrt{k}}\right)^2=D_X\left(\frac{w}{\sqrt{k}},q,\theta
\right)^2 =D_L(q)^2+O(k^{-1/2})$, where $\mathrm{dens}_L=D_L^2\cdot
|dq|$. Since $\int
_{\mathbb{R}^\mathrm{n}}e^{-2\|p\|^2}\,dp=(\pi/2)^{\mathrm{n}/2}$,
given (\ref{eqn:asympt-exp-product-diff}) we obtain
 $A_{jk}\sim
 k^{2+\mathrm{n}/2}\cdot r^2\left (\frac{2}{\pi}\right
 )^{\mathrm{n}/2}\int _L\varphi _j\,F\big(\widetilde{W},\widetilde{W}'\big)\cdot
 \mathrm{dens}_L+\sum _{h\ge 1}u_{jh}k^{2+(\mathrm{n}-h)/2}$.
Summing over $j$, we get
(\ref{eqn:asympt-exp-hermitian-product-phi-k}). \hfill $\square$

\section{Proof of Theorem \ref{thm:semicl-asympt-isom}.}

Suppose that $\mathfrak{I}\subseteq \mathbf{L}_{\mathrm{BS}}$ is the
isodrastic leaf such that
$\mathfrak{S}=\mathbf{W}_{\mathrm{h}}\mathfrak{I}$. Given
$(L,\lambda)\in \mathbf{W}_{\mathrm{h}}\mathfrak{I}$, we shall now
consider the asymptotics of the derivative of the projectivized BPU
map,
$d_{(L,\lambda)}\Phi_k:T_{(L,\lambda)}
\mathbf{W}_{\mathrm{h}}\mathfrak{I}\rightarrow T_{\Phi
_k(L,\lambda)}\mathbb{P}H(X)_k$,
for $k=l\,r$ and $l\rightarrow +\infty$.

Let $\big(V,\left<\,\,\right>\big)$ be a finite-dimensional unitary
vector space. Let $\varsigma:V\setminus \{0\}\rightarrow
\mathbb{P}V$, $v\mapsto [v]$ be the projection. If $v\in V\setminus
\{0\}$, the differential $d_v\varsigma:V\rightarrow
T_{[v]}\mathbb{P}V$ induces an algebraic isomorphism
$v^\perp\rightarrow T_{[v]}\mathbb{P}V$, where $v^\perp\subseteq V$
is the unitary orthocomplement of $v$. The Fubini-Study metric is
determined by
\begin{equation}\label{eqn:unitary-product-on-projective}
\left <d_v\varsigma(v'),d_v\varsigma(v'')\right >_{[v]}= \frac{\left
<v',v''\right >}{\left <v,v\right >}\,\,\,\,\,\,\,\,\,(v',v''\in
v^\perp).\end{equation}

Suppose $W=W\big(f,\ell\big)\in
T_{(L,\lambda)}\mathbf{W}_{\mathrm{h}}\mathfrak{I}$, where
$\int _Lf\,\lambda\bullet\lambda=\int _L\ell \bullet \lambda=0$,
and set $W_k=:W(f,k\ell)$. Since by assumption $L$ is a BSL
submanifold, by definition $\exists\,P\subseteq X$ compact and
connected Planckian submanifold  with $L=\pi(P)$. Thus
$(P,\lambda)\in \mathbf{W}_{\mathrm{h}}\mathbf{P}_\mathfrak{I}$ lies
over $(L,\lambda)$.
Now $W_k$ naturally lifts to $\widetilde{W}_k\in T_{(P,\lambda)}
\mathbf{W}_{\mathrm{h}}\mathbf{P}_\mathfrak{I}$, hence
$W_k=d_{(P,\lambda)}\widehat{\pi}\big (\widetilde{W}_k\big)$, where
$\widehat{\pi}:\mathbf{W}_{\mathrm{h}}\mathbf{P}_\mathfrak{I}\rightarrow
\mathbf{W}_{\mathrm{h}}\mathfrak{I}$ is the projection. If
$\varsigma:H(X)_k\setminus \{0\}\rightarrow \mathbb{P}H(X)_k$ is the
projection, then $\Phi _k\circ \widehat{\pi}=\varsigma \circ
\widetilde{\Phi}_k$. Therefore,
\begin{equation}\label{eqn:diff-linear-non-linear-bis}
 d_{(L,\lambda)}\Phi_k(W_k)=d_{(L,\lambda)}\Phi_k\circ d_{(P,\lambda)}
 \widehat{\pi}\big (\widetilde{W}_k\big)
 =d_{\widetilde{\Phi}_k(P,\lambda)}\varsigma\Big (d_{(P,\lambda)}
 \widetilde{\Phi}_k(\widetilde{W}_k)\Big).
 \end{equation}
Now $\left
(\widetilde{\Phi}_k(P,\lambda),d_{(P,\lambda)}\widetilde{\Phi}_k(\widetilde{W_k})
\right)_{L^2(X)} \sim \sum _{h\ge 0}s_h\,k^{(1+\mathrm{n}-h)/2}$, in
view of (\ref{eqn:asympt-exp-hermitian-product-phi-k-dphi-k}) and
the conditions on $f$ and $\ell$.
Define $Z_k\in H(X)_k$, $k\gg 0$, by
\begin{equation}\label{eqn:z-k}
Z_k=:d_{(P,\lambda)}\widetilde{\Phi}_k(\widetilde{W}_k)- \left
[\frac{\left(\widetilde{\Phi}_k(P,\lambda),d_{(P,\lambda)}\widetilde{\Phi}_k\big(\widetilde{W}_k\big)\right)
_{L^2(X)}} {\left (\widetilde{\Phi}_k(P,\lambda),
\widetilde{\Phi}_k(P,\lambda)\right)_{L^2(X)}} \right]\cdot
\widetilde{\Phi}_k(P,\lambda).\end{equation}
Thus, $\big(Z_k,\widetilde{\Phi}_k(P,\lambda)\big)_{L^2(X)}=0$, and
furthermore
$d_{(L,\lambda)}\Phi_k(W_k)=d_{\widetilde{\Phi}_k(P,\lambda)}\varsigma(Z_k)$
by (\ref{eqn:diff-linear-non-linear-bis}). Suppose now that
$W'=W(f',\ell')\in
T_{(L,\lambda)}\mathbf{W}_{\mathrm{h}}\mathfrak{I}$ is a second
tangent vector, and let $W_k'$ and $Z_k'$ be defined as $W_k$ and
$Z_k$, starting from $W'$. Then using Proposition
\ref{prop:asymptotically-orthogonal} and the above we obtain an
asymptotic expansion
\begin{equation*}
\label{eqn:Z-k-paired-Z-prime-k} \left (Z_k,Z'_k\right)_{L^2(X)}\sim
k^{2+\mathrm{n}/2}\cdot r^2\left (\frac{2}{\pi}\right
 )^{\mathrm{n}/2}\int _LF\big(\widetilde{W},\widetilde{W}'\big)\cdot
 \mathrm{dens}_L+\sum _{h\ge 1}L_{h}k^{2+(\mathrm{n}-h)/2}.
\end{equation*}
Again in view of Proposition \ref{prop:asymptotically-orthogonal},
we deduce from (\ref{eqn:unitary-product-on-projective}):
\begin{eqnarray}
\label{eqn:asympt-exp-for-differential-and-pull-back} \lefteqn{
\Big(d_{(L,\lambda)}\Phi_k(W_k),d_{(L,\lambda)}\Phi_k(W_k')\Big)_{\Phi_k(L,\lambda)}
=\left (d_{\widetilde{\Phi}_k(P,\lambda)}\varsigma(Z_k),
d_{\widetilde{\Phi}_k(P,\lambda)}\varsigma(Z_k')\right)
_{[\widetilde{\Phi}_k(P,\lambda)]} }
\\
&&=\frac{\left (Z_k,Z_k'\right)_{L^2(X)}}{\left
(\widetilde{\Phi}_k(P,\lambda),\widetilde{\Phi}_k(P,\lambda)\right)_{L^2(X)}}
\sim k^{2}\cdot \int _LF\big(\widetilde{W},\widetilde{W}'\big)\cdot
 \mathrm{dens}_L
+\sum _{h\ge 1}L_{h}'k^{2-h/2}.\nonumber
\end{eqnarray}
Given Remark \ref{rem:integral-of-F-is-product}, to complete the
proof of Theorem \ref{thm:semicl-asympt-isom} we need only take real
and imaginary parts in
(\ref{eqn:asympt-exp-for-differential-and-pull-back}).
\hfill
$\square$

\end{document}